\documentclass[12pt,onesite]{article}
\usepackage{amsmath,amsxtra,amssymb,latexsym,amscd}

\numberwithin{equation}{section}

\font\abst=cmr10

\newtheorem{The}{Theorem}[section]

\newtheorem{Cor}{Corollary}[section]

\date{}

\begin{document}
\begin{center}
\Large{\bf  A note on the combination of equilibrium problems}\\
\vspace*{1cm}
{\small \sc{Nguyen Thi Thanh Ha}\footnote{Email: nttha711@gmail.com}},
{\small \sc{Tran Thi Huyen Thanh}\footnote{Email: thanh0712@gmail.com}}, \\
{\small \sc{Nguyen Ngoc Hai}\footnote{Email: hainn@dhcd.edu.vn},}
{\small \sc{Hy Duc Manh}\footnote{Email:manhhd@gmail.com}},
{\small \sc{Bui Van Dinh}\footnote{Corresponding Email: vandinhb@gmail.com\vspace*{0.1cm}}\\
\small{ $^{1,2,4,5}${\it Department of Mathematics, Le Quy Don Technical University, Hanoi, Vietnam}}\\
\small{\noindent $^3${\it Department of Scientific Fundamentals, Trade Union University, Hanoi, Vietnam}}}
 \end{center}

\noindent {\bf Abstract.}
\noindent In this short paper, we show that the solution set of a combination of equilibrium problems is not necessary contained in the intersection of a finite family of solution sets of equilibrium problems. As a corollary, we deduce that statements in recent papers given by S. Suwannaut, A. Kangtunyakarn (Fixed Point Theory Appl. 2013, 2014; Thai Journal of Maths. 2016), W. Khuangsatung, A. Kangtunyakarn (Fixed Point Theory Appl. 2014),  and  A.A. Khan, W. Cholamjiak, and K.R. Kazmi (Comput. Appl. Maths. 2018) are not correct.


\noindent{\bf Keywords and Phrase:}  {\abst equilibria; Ky Fan inequality; Combination.}

\section{Introduction}
\indent Let $C$ be a nonempty closed convex subset in the Euclidean space $\mathbb{R}^n$ and 
 $f : C \times C \to \mathbb{R}$ be a bifunction. The equilibrium problem (shortly EP($C, f$)), in the sense of Blum, Muu and Oettli \cite{BO, MO} (see also \cite{Fan}), consists of finding $ x^* \in C $ such that
\begin{equation}\notag
  f(x^*, y) \geq 0, \  \forall y \in C.
\end{equation}
 We denote the solution set of EP($C, f$) by $Sol(C, f)$. Solution methods for EP($C,f$) can be found in \cite{QMH,DK}.


Let $f_i : C \times C \to \mathbb{R}, i = 1, 2, ..., N$,  be bifunctions defined on $C$. Recently, many researchers are interested in finding a common solution of a finite family of equilibrium problems \cite{SK13,SK14,SK16,KCK} (CSEP for short).
\begin{equation}\notag
  \text{ Find } x^* \in C \text{ such that } f_i(x^*, y) \geq 0, \  \forall y \in C \text{ and } i = 1, 2, ..., N. \eqno{CSEP(C, f_i)}
\end{equation}
Or, equivalently, $$   \text{ find } x^* \in \Omega := \cap_{i=1}^N Sol(C,f_i).
 $$
Given bifunctions $f_i$, $i = 1,..., N$ defined on $C$. Let $\alpha_i \in (0, 1), i = 1, ..., N$ such that $\sum_{i = 1}^N \alpha_i = 1.$ Set
$$ f(x, y) = \sum_{i = 1}^N \alpha_i f_i(x, y).$$
The combination of equilibrium problems (shortly, CEP$(C, \sum_{i = 1}^N \alpha_i f_i)$) consists of finding $x^* \in C$ such that
$$ f(x^*, y) = \sum_{i = 1}^N \alpha_i f_i(x^*, y) \geq 0, \forall y \in C.$$
By $Sol(C, \sum_{i = 1}^N \alpha_i f_i)$, we denote the solution set of the combination of equilibrium problems. \\
 In 2013, S. Suwannaut and A. Kangtunyakarn \cite{SK13} said that under certain conditions then
$$ \Omega := \cap_{i=1}^N Sol(C,f_i) =  Sol(C, \sum_{i = 1}^N \alpha_i f_i). $$
  Therefore, to find  a common solution of a finite family of equilibrium problems leads to find a solution of a combination of equilibrium problems CEP($C, \sum_{i = 1}^N \alpha_i f_i$).
Based on this relation, S. Suwannaut and Kangtunyakarn \cite{SK13,SK14,SK16}, W. Khuangsatung and A. Kangtunyakarn \cite{KK},  S.A. Khan, W. Cholamjiak, and K.R. Kazmi \cite{KCK} gave algorithms for finding a common element of the fixed point sets of a family of mappings and the solution sets of equilibrium problems and/or the zero point sets of a family of mappings.

In this short paper, we show that, under the same conditions given in \cite{SK13}, the relation 
$$ Sol(C, \sum_{i = 1}^N \alpha_i f_i) \subset \cap_{i=1}^N Sol(C,f_i), $$
does not hold true. Therefore, presenting of recent papers \cite{SK13,SK14,SK16,KK,KCK} using this formula are not correct.

The rest of paper is organized as follows. The next section contains some preliminaries on  equilibrium problems and some statements in papers \cite{SK13,SK14,SK16,KK,KCK}  related with combination of equilibrium problems. The last section is devoted to show that the common points of a finite family of equilibrium problems is truly contained in a solution set of a combination of equilibrium problems and its corollaries.


\section{Preliminaries}
In this section, we present some statements presented in recent papers related to combination of equilibrium problems. Let $\varphi : C \times C \to \mathbb{R}$ be a bifunction defined on $C$. In the sequel, we need the following blanket assumptions:

\noindent{\bf Assumptions $\mathcal{A}.$}

\begin{itemize}
\item[($\mathcal{A}_1$)]  $\varphi(x, x) = 0$ for  every $x\in C$;
\item[($\mathcal{A}_2$)]  $\varphi$ is monotone on $C$;
\item[($\mathcal{A}_3$)] $\varphi$ is upper hemicontinuous, i.e., for each $x, y, z \in C$ we have
$$\lim\sup_{t \to 0^+} \varphi (tz + (1-t)x, y) \leq \varphi (x,y);$$
\item[($\mathcal{A}_4$)] for each $x \in C$, $\varphi (x, \cdot)$ is lower semicontinuos and convex on $C$;
\item[($\mathcal{A}_5$)] for fixed $r > 0$ and $z \in C$, there exists a nonempty compact convex subset $B$ of $\mathbb{R}^n$ and $x \in C \cap B$, such that
$$\varphi(y,x) + \frac{1}{r} \langle y-z, z-x \rangle < 0, \forall y \in C \setminus\ B.$$
\end{itemize}

The following statement is in \cite{SK13}.

\noindent{\bf Statement 2.1} (See \cite[Lemma 2.7]{SK13}).  
Let  $f_i, i = 1, 2, ..., N$ be bifunctions satisfying $\mathcal{A}_1 - \mathcal{A}_4$ with $\cap_{i=1}^N Sol(C, f_i) \neq \emptyset$. Then
$$\cap_{i=1}^N Sol(C, f_i) = Sol(C, f),$$
where $f(x,y) = \sum_{i=1}^N \alpha_i f_i(x,y)$, $\alpha_i >0, \forall i = 1, 2, ..., N$ and $\sum_{i=1}^N \alpha_i = 1$.

If Statement 2.1 holds true then it allows us to find  common  solutions of  $N$ equilibrium problems by solving  a combination of equilibrium problems.


The following statement is in \cite{SK14}. \\

\noindent{\bf Statement 2.2} (See \cite[Theorem 3.1]{SK14}). Let $F$ be an an $\tau$-contractive mapping on $\mathbb{R}^n$ and let $A$ be a strongly positive linear bounded operator on $\mathbb{R}^n$ with coefficient $\bar\gamma$ and $0 < \gamma < \frac{\bar\gamma}{\tau}$. For every $i = 1, 2, ..., N$ let $f_i : C \times C \to \mathbb{R}$ be a bifunction satisfying $\mathcal{A}_1 - \mathcal{A}_4$ with $\Omega=\cap_{i=1}^N Sol(C, f_i) \not=\emptyset$. Let $\{x^n\}, \{y^n\}, \{z^n\}$ be sequences generated by $x^1 \in \mathbb{R}^n$ and
\begin{equation*}
\begin{cases}
\sum_{i=1}^N \alpha_i f_i(z^k,y) + \frac{1}{\rho_k}\langle y-z^k, x^k-z^k\rangle \geq 0, \forall y \in C, \\
y^k = \theta_k P_C(x^k) + (1-\theta_k)z^k,\\
x^{k+1} = \delta_k\gamma F(x^k)+(I-\delta_kA)y^k, 
\end{cases}
\end{equation*}
where $\{\delta_k\}, \{\theta_k\}, \{\rho_k\} \subset (0, 1), 0< \alpha_i < 1, \forall i = 1, ..., N$. Suppose the conditions $(i)-(vi)$ hold.
\begin{itemize}
\item[(i)] $\lim_{k \to \infty} \delta_k = 0$ and $\sum_{k=0}^{\infty} \delta_k = \infty$;
\item[(ii)] $0 < \underline{\theta} \leq \theta_k \leq \bar{\theta} < 1$, for some $\underline{\theta}, \bar{\theta} \in (0, 1)$;
\item[(iii)] $0 < \underline{\alpha} \leq \alpha_k \leq \bar{\alpha} < 1$, for some $\underline{\alpha}, \bar{\alpha} \in (0, 1)$;
\item[(iv)] $\sum_{i=1}^N \alpha_i = 1$; 
\item[(v)]$\sum_{i=1}^N |\delta_{k+1} - \delta_k| < \infty$, $\sum_{i=1}^\infty |\theta_{k+1} - \delta_k| < \infty$, $\sum_{i=1}^\infty |\rho_{k+1} - \rho_k| < \infty$. 
\end{itemize}
Then the sequences $\{x^k\}, \{y^k\}$, and $\{z^k\}$ converge to $q = P_\Omega(I-A+\gamma F)q$.

From Theorem 3.1 in \cite{KK} we get the following statement.

\noindent{\bf Statement 2.3} (See \cite[Theorem 3.1]{KK}). Let $f_i, i = 1, 2, ..., N$ satisfy assumption $\mathcal{A}_1 - \mathcal{A}_4$.
Assume that $\Omega =  \cap_{i=1}^N Sol(C, f_i) \neq \emptyset$. Let the sequence $\{x^k\}$ and $\{y^k\}$ be generated by $u, x^1 \in \mathbb{R}^n$ and
    \begin{equation*}
\begin{cases}
\sum_{i=1}^N \alpha_i f_i(y^k,y) + \frac{1}{\rho_k}\langle y-y^k, y^k-x^k\rangle \geq 0, \forall y \in C, \\
x^{k+1} = \lambda_k u +\mu_k x^k +\delta_k y^k 
\end{cases}
\end{equation*}
where $\{\lambda_k\}, \{\mu_k\}, \{\delta_k\} \subset (0, 1)$ and $\lambda_k +\mu_k+ \delta_k =1$; $\{\rho_k\} \subset (\underline{\rho}, \bar{\rho}) \subset (0, 1)$, $ 0< \alpha_i < 1, \forall i = 1, ..., N$. Suppose the conditions $(i)-(iii)$ hold.
\begin{itemize}
\item[(i)] $\lim_{k \to \infty} \lambda_k = 0$ and $\sum_{k=0}^{\infty} \lambda_k = \infty$;
\item[(ii)] $\sum_{i=1}^N \alpha_i = 1$; 
\item[(iii)]$\sum_{i=1}^N |\delta_{k+1} - \delta_k| < \infty$.
\end{itemize}
Then the sequences $\{x^k\}, \{y^k\}$ converge to $q = P_\Omega (u)$.

The next statement is deduced from Theorem 3.1 in \cite{SK16}.

\noindent{\bf Statement 2.4} \cite[Theorem 3.1]{SK16}. Let $F$ be an an $\tau$-contractive mapping on $\mathbb{R}^n$ and let $f_i, i = 1, 2, ..., N$ satisfy assumption $\mathcal{A}_1 - \mathcal{A}_4$.
Assume that $\Omega =  \cap_{i=1}^N Sol(C, f_i) \neq \emptyset$. Let the sequence $\{x^k\}$ and $\{y^k\}$ be generated by $ x^1 \in C$ and
    \begin{equation*}
\begin{cases}
\sum_{i=1}^N \alpha_i f_i(y^k,y) + \frac{1}{\rho_k}\langle y-y^k, y^k-x^k\rangle \geq 0, \forall y \in C, \\
x^{k+1} = \lambda_k F(x^k) +\mu_k x^k +\delta_k y^k 
\end{cases}
\end{equation*}
where $\{\lambda_k\}, \{\mu_k\}, \{\delta_k\} \subset (0, 1)$; $\{\rho_k\} \subset (\underline{\rho}, \bar{\rho}) \subset (0, 1)$, $ 0< \alpha_i < 1, \forall i = 1, ..., N$. Suppose the conditions $(i)-(iii)$ hold.
\begin{itemize}
\item[(i)] $\lim_{k \to \infty} \lambda_k = 0$ and $\sum_{k=0}^{\infty} \lambda_k = \infty$;
\item[(ii)] $\sum_{i=1}^N \alpha_i = 1$; 
\item[(iii)]$\sum_{i=1}^{\infty} |\rho_{k+1} - \rho_k| < \infty$.
\end{itemize}
Then the sequences $\{x^k\}, \{y^k\}$ converge to $q = P_\Omega (u)$.

From Theorem 4.2 in \cite{KCK} we get the following statement.

\noindent{\bf Statement 2.5} \cite[Theorem 3.1]{KCK}. 
Let $f_i, i = 1, 2, ..., N$ satisfy assumption $\mathcal{A}$.
Assume that $\Omega =  \cap_{i=1}^N Sol(C, f_i) \neq \emptyset$. For given $x^0, x^1 \in \mathbb{R}^n$, let the sequence $\{x^k\}$,  $\{y^k\}$ and $z^k$ be generated by 
    \begin{equation*}
\begin{cases}
y^k = x^k + \theta_k (x^k - x^{k-1})\\
\sum_{i=1}^N \alpha_i f_i(z^k,y) + \frac{1}{\rho_k}\langle y-z^k, z^k-y^k\rangle \geq 0, \forall y \in C, \\
x^{k+1} = \lambda_k x^k +\mu_k z^k 
\end{cases}
\end{equation*}
where $\{\theta_k\} \subset [0, \theta], \theta \in [0; 1]$, $\{\lambda_k\}, \{\mu_k\} \subset (0, 1)$ and $\lambda_k + \mu_k = 1$ for all $k$; $\{\rho_k\} \subset (\underline{\rho}, \bar{\rho}) \subset (0, 1)$, $ 0< \alpha_i < 1, \forall i = 1, ..., N$. Suppose that the following conditions hold
\begin{itemize}
\item[(i)] $\theta_k \|x^k - x_{k-1}\| < \infty$;
\item[(ii)] $\sum_{i=1}^{\infty} \alpha_i < \infty$ and $\lim_{i \to \infty} \alpha_i = 0$; 
\item[(iii)] $\sum_{i=1}^{\infty} |\rho_{k+1} - \rho_k| < \infty$, $\sum_{i=1}^{\infty} |\lambda_{k+1} - \lambda_k| < \infty$.
\end{itemize}
Then the sequence $\{x^k\}$ converges to $q = P_\Omega (u)$.


\section{Main Results}
Now, given natural number $N \geq 2$ and a nonempty, closed convex set $C$ and bifunctions $f_i$ ($i = 1...N$) defined on $C$ such that $$\cap_{i=1}^N Sol(C, f_i) \neq \emptyset.$$
 For $\alpha_i \in (0, 1), i = 1, ..., N$ and $\sum_{i=1}^N \alpha_i = 1$. We define
 $$ f(x,y) = \sum_{i=1}^N \alpha_i f_i(x,y).$$
It is clear that if $x^* \in \cap_{i=1}^N Sol(C, f_i)$ then $f_i(x^*, y)  \geq 0, \forall y \in C, i = 1, 2, ..., N.$ Therefore $f(x^*, y) = \sum_{i=1}^N \alpha_i f_i(x^*,y)\geq 0, \forall y \in C$. So $x^* \in Sol(C, f)$.\\
Hence 
$$\cap_{i=1}^N Sol(C, f_i) \subset Sol(C,f).$$

 The following theorem show that under assumptions $\mathcal{A}_1-\mathcal{A}_4$, the inversion is not true.
 
\begin{The}\label{T4.1} For any integer number $N \geq 2$, there exist a nonempty, closed convex set $C$ and bifunctions $f_1, f_2, ..., f_N$ defined on $C$ satisfy assumptions $\mathcal{A}_1-\mathcal{A}_4$  and $\alpha_i \in (0, 1), i= 1, 2, ..., N, \sum_{i=1}^N \alpha_i = 1$ such that
$$Sol(C, \sum_{i=1}^N \alpha_i f_i) \not\subset \cap_{i=1}^N Sol(C, f_i).$$
  \end{The}
{\it Proof.} 
It is clear that, we only need prove for the case $n=2$ and $N=2$. Indeed, for  $x = (x_1, x_2)\in \mathbb{R}^2$, $y =(y_1, y_2) \in \mathbb{R}^2$. Consider the set $C$ and bifunctions are given as follow
 $$C = \{ (x_1, x_2) \in \mathbb{R}^2 : x_1 \geq 0, x_2 \geq 0\}.$$
$$f_1(x,y) = x_2y_1 - x_1y_2.$$ 
$$f_2(x,y) = x_1y_2 - x_2y_1.$$ 
Then we have: $f_1(x,x) = 0, \forall x \in C$. For all $x, y \in C$, we have
$$f_1(x,y)+ f_1(y,x) = x_2y_1 - x_1y_2 + y_2x_1 - y_1x_2 = 0.$$
Hence, $f_1$ is monotone on $C$.

For each $x \in C, f(x, y)$ is linear in $y$, so $f(x, \cdot)$ is convex. It is trivial that $f_1$ is continuous on $C \times C$.

 Therefore bifunction $f_1$ satisfies assumptions $\mathcal{A}_1, \mathcal{A}_2, \mathcal{A}_3,$ and $ \mathcal{A}_4$.
 
Similarly, $f_2$ satisfies assumptions $\mathcal{A}_1, \mathcal{A}_2, \mathcal{A}_3,$ and $ \mathcal{A}_4$.
 In addition, It can be seen that
 $$Sol(C, f_1) = \{0\} \times [0, +\infty).$$
 $$Sol(C, f_2) =  [0, +\infty) \times \{0\}.$$
So, $$Sol(C, f_1) \cap Sol(C, f_2) = \{(0,0)\}.$$
Now, we consider a combination of $f_1, f_2$ given as follows
\begin{equation*}
\begin{aligned}
f(x,y) &= \frac{1}{2} f_1(x,y) + \frac{1}{2} f_2(x,y)\\
       & = \frac{1}{2} \big[ f_1(x,y) + f_2(x,y)\big]\\
       & = 0, \forall x, y \in C.
\end{aligned}
\end{equation*}
It is obvious that $f$ satisfies assumptions $\mathcal{A}_1, \mathcal{A}_2, \mathcal{A}_3,$ and $ \mathcal{A}_4$. Moreover
$$Sol(C, f) = C = [0, +\infty) \times [0, +\infty).$$
Therefore $$Sol(C,f) \not\subset Sol(C, f_1) \cap Sol(C, f_2).$$
\hfill$\Box$

From this theorem, we have the following corollary
\begin{Cor}
 Statement 2.1 - Statement 2.5 are not correct.
\end{Cor}
{\bf Proof.} We take $N=2$, the set $C$, bifunctions $f_1$ and $f_2$ defined as in Theorem~\ref{T4.1}. The combination of $f_1$ and $f_2$ is given by $f(x,y) = \frac{1}{2} f_1(x,y) + \frac{1}{2} f_2(x,y) = 0, \forall x, y \in C$. Hence, $\Omega = Sol(C, f_1) \cap Sol(C,f_2) = \{(0, 0)\}$,  $Sol(C, f) = C = [0, +\infty) \times [0, +\infty).$   Then we have the followings:
\begin{itemize}
\item[(a)] Statement 2.1 is false.
\item[(b)]  Take $x^1 \in C$ such that $x^1 \neq (0,0)$ and set $F(x) = x^1, Ax = x, \  \forall x \in \mathbb{R}^n$. Choose $\gamma =1$, then the sequence $\{x^k\}$ generated by Statement 2.2 takes the form
$$x^{k+1} = \delta_k x^1 + (1 - \delta_k)x^k = x^1, \forall k.$$
Therefore, it converses to $x^1 \not\in \Omega.$ It means that Statement 2 is false.
\item[(c)] By taking any $u= x^1 \in C$ such that $x^1 \neq (0,0)$. Then the sequence $\{x^k\}$ generated by the scheme in Statement 2.3 becomes
$$x^{k+1} = \lambda_k x^1 + (1-\lambda_k)x^k = x^1, \forall k.$$
It leads to $x^k \to x^1 \not\in \Omega$. Hence Statement 2.3 is not correct.
\item[(d)] Similar to the case (b), we have Statement 2.4 is false.
\item[(e)] By taking any $x^1 = x^0 \in C$, then the sequence $\{x^k\}$ generated by Statement 2.5 takes the form
$$x^{k} = x^1, \forall k.$$
So, Statement 2.5 does not true.
\end{itemize}

\hfill$\Box$


\newpage

 \noindent{\bf Conclusion}. We have proved that there exist a finite family of monotone equilibrium problems such that the common solution set of them does not contain the solution set of a combination of those equilibrium problems. Based on this fact, we imply that recent papers \cite{KCK,KK,SK13,SK14,SK16} are not correct.


\bigskip
\bibliographystyle{amsplain}

\end{document}